# The decolonisation of mathematics

John Armstrong and India Jackman

**Abstract:** *We describe a mainstream "universalist" approach to the understanding of mathematics. We then conduct a systematic (but not exhaustive) review of the academic literature on the decolonisation of mathematics and identify how this challenges the universalist view. We examine evidence of whether the experience of mathematics in the UK is systemically racist, examining both the decolonial arguments and the empirical evidence. We find that there may be some benefit in teaching the history of mathematics, but that this should be weighed against the opportunity cost. We find some prima-facie evidence of discrimination in the descriptive statistics on the representation of ethnic minorities in academic roles in UK higher education.*

## Introduction

In the UK, the movement to decolonise the curriculum has gained considerable momentum in recent years, driven by the Black Lives Matter movement. Decolonisation of the curriculum had its origins in the humanities, but writing a guest editorial in Nature, Nobles, Womack, Wonkam, and Wathuti (2022) state "It is so important for science curricula, research and academic spaces to go through decolonization processes. These are not political or ideological acts, but part of science itself — an example of science's self-correcting mechanism in the pursuit of truth."

In popular understanding, decolonising the curriculum means ensuring that the reading lists include works by black and minority authors, and that such authors are not treated as being inferior to white male authors. But it is less obvious what decolonisation means in relation to mathematics. This is acknowledged in another editorial in Nature (2023), "One common argument is that decolonization is irrelevant to the practice of mathematics: the solution to a quadratic equation doesn't, after all, depend on a mathematician's identity or protected characteristics". They go on to state that to understand decolonisation we must first "reprise aspects of an older, more academically focused debate on whether — or to what extent — scientific knowledge is socially constructed".

With this in mind, our paper begins with a review of a conventional mathematical understanding of the history and philosophy of mathematics. If we are to teach mathematics from a decolonial perspective, it is important that we can explain both the decolonial and conventional perspective to students, and how they differ. Therefore, this section has been written for a reader with no specialist experience in mathematics.

The second part of the paper presents the challenges to this conventional viewpoint made by the decolonisation agenda. We review the key texts on the decolonisation of mathematics grouped into the following themes: decolonial theory; critical pedagogy; indigenous minorities; Africa; and the UK. To avoid giving a skewed presentation, we augment this discussion of key texts with a review of a sample of 37 texts taken systematically from the literature on the decolonisation of mathematics. Having reviewed the literature, we will argue that decolonisation should not be seen simply a matter of promoting the education of black and minority students; rather, it requires challenging conventional understandings of mathematics. We will present evidence that this is a political and ideological act rather than a part of science.

It is possible to believe that mathematics as a discipline is systemically racist whether or not one supports the decolonisation of mathematics. In the third part of the paper, we examine the evidence for systemic racism in mathematics from a UK viewpoint. We will examine the decolonial argument that mathematics practice is racist, but will also consider other forms of evidence including data on minority presentation among university academics, the contents of GCSE textbooks and a Masters level module taught by one of the authors.

We conclude with a discussion of the overall strength of the evidence base for decolonising the mathematics curriculum in the UK and we assess the desirability of this programme.

To the best of our knowledge, this paper is the first to review the literature on the decolonisation of mathematics. We hope that our work will inform the debate not only on the decolonisation of mathematics, but also on the wider programme of decolonisation of the curriculum. For example, the "indigenization" of science is debated in a series of essays in (Widdowson, 2016) but, as Widdowson acknowledges that it was difficult to get supporters of indigenization to contribute to this volume, one can't assume that the views expressed are representative. However, our review of the mathematical literature identifies similar themes to those discussed in Widdowson, suggesting that it succeeds in giving a fair representation of decolonial perspectives.

Our central aim in this paper is to show the importance of the question of whether mathematical knowledge is socially constructed to the decolonisation of mathematics. We do not claim to break new ground on answering the question of whether mathematics is socially constructed, preferring as far as possible to provide the reader with the tools needed to decide for themselves. We focus on mathematics as the most challenging test-case for the decolonisation programme. Authors such as (Brown, 2009; Sokal & Bricmont, 1998) present arguments against the social construction of science more broadly, but this is beyond the scope of this paper.

## Part 1: A universalist account of mathematics

Before discussing decolonisation, we would like to give an account of some key aspects of the history, philosophy and practice of mathematics that assumes a "universalist" perspective. We will see how this view is challenged in the second part of the paper. Where we do not provide explicit references to the history one should refer to a modern history of mathematics such as (Boyer & Merzbach, 2011).

The number system we use today was developed between the first and fourth centuries AD by Indian mathematicians and became known in Europe via the works of the Persian mathematician Al-Khwareizmi. Before this, the Babylonian civilization had developed a base-60 number system and the Chinese had developed a base-10 system of counting rods. Independently, in the pre-Colombian Americas, the Maya used a base -20 counting system, possibly developed by the Olmecs. Mathematics has proved extraordinarily successful and has been used to facilitate trade and administration in Ancient Egypt, Tang dynasty China, pre-Colombian America, democracies, monarchies, communist states, dictatorships, empires, theocracies and republics.

Mathematical theory has also developed internationally. In 14[th]-15[th] century Kerala, Mādhava of Sangamagrāma, or perhaps one of his followers, proved the following infinite series formula for $\pi$:

$$\frac{\pi}{4} = \frac{1}{1} - \frac{1}{3} + \frac{1}{5} - \frac{1}{7} + \frac{1}{9} - \frac{1}{11} + \cdots$$

This formula was then derived apparently independently in Europe, first by Gregory in 1671, then Leibniz in 1673. While, it cannot be ruled out that Gregory or Leibniz were influenced by the Keralan school, despite searches no direct evidence has been found of this influence (Joseph, 2010; Roy, 1990).

Mādhava's formula gives the value of π, a fundamentally geometric concept which can be defined as the ratio of a circle's circumference to its diameter. The right-hand side of the formula is arithmetic rather than geometric: we have the sequence of odd numbers and the arithmetic operations of addition, subtraction and division. It seems that Mādhava and Leibniz shared common conceptions of geometry and arithmetic despite not sharing either language or culture: Mādhava used very different notation from Leibniz, yet Mādhava and Leibniz would have agreed that their results were essentially equivalent.

One partial explanation of these common conceptions is that the laws of physics do not differ markedly between Germany and Kerala, hence similar results will be found about the properties of circles wherever they are studied. However, this is only a partial explanation because the Mādhava-Leibniz formula is not a result about physical circles at all. Their formula involves an infinite sequence of numbers and is intended to indicate the value of π exactly. Yet since it is physically impossible to measure the circumference of a circle with arbitrary accuracy, it surely does not make sense to state the circumference of any physical circle with infinite precision.

If a mathematical circle is not a physical circle, what exactly is geometry? The Ancient Greek philosopher Plato argued that geometry was the study of the non-physical underlying essence of objects such as a circle which he termed their Ideal: the physical objects we see are an imperfect shadow of this underlying Ideal. Plato's belief is highly metaphysical and it is not clear what is gained in practice from this theory, and many mathematicians reject Platonism or see it as irrelevant.

One Ancient Greek idea that does enjoy broad support among mathematicians is the idea of mathematical proof. Even if we are not sure what the metaphysical meaning of "circle" is, if we are willing to make certain assumptions about a circle's properties then we may be able to say what logical consequences flow from these assumptions.

Not all readers will have seen a mathematical proof, so we give as an example a proof of Pythagoras's Theorem.

**Pythagoras's Theorem:**

**Claim:** If one draws squares on the edges of a right-angle triangle as shown in Figure 1, then the area of the large square is equal to the sum of the two smaller squares.

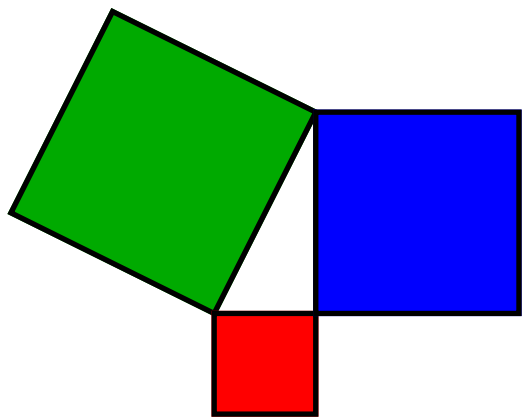

**Figure 1: Arrangement of squares in Pythagoras's theorem.**

**Proof:**

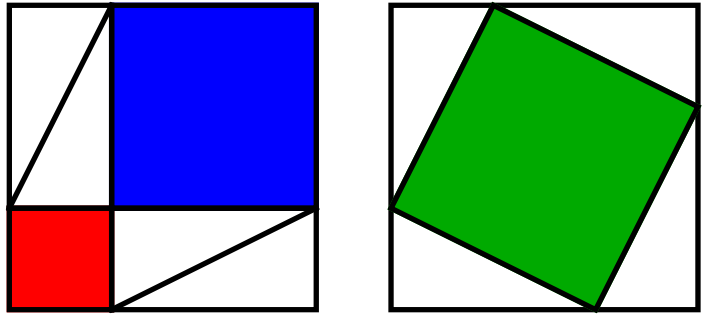

**Figures 2a, 2b: A geometric proof of Pythagoras's theorem by rearrangement.**

Draw 4 copies of the triangle and place them as shown in white in Figure 2a. Draw another 4 copies and place them as shown in white in Figure 2b. Figure 2a and Figure 2b are both pictures of squares of the same size, so they have the same area. Since the same white triangles appear in Figure 2a and Figure 2b, this means that the areas of the green square in Figure 2b is equal to the sum of the areas of the coloured squares in Figure 2a. Matching the areas in Figures 2a and 2b to the corresponding areas in Figure 1 gives the result.

The elegance of this rearrangement proof is an example of what mathematicians mean by "beauty". The history of Pythagoras's theorem is debated and it was most likely known by the Babylonians more than 1000 years before Pythagoras (Neugebauer, 1969). A statement of the theorem is given in the Indian Baudhayana Shulba Sutra (c. 100 BC). A pictorial proof of one special case is shown in Figure 3, taken from the Chinese Zhoubi Shuanjing (unknown date, perhaps c. 100 BC).

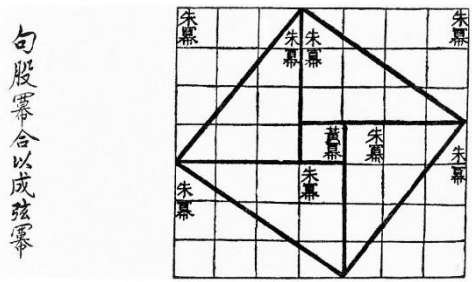

**Figure 3:** A visual proof of a special case of Pythagoras's Theorem taken from the Zhoubi Shuanjing

The rearrangement proof makes many assumptions. It assumes that the pictures we have drawn are accurate, and that the shapes we have drawn that look like squares really are squares. It assumes that if you divide a shape into a number of smaller pieces whose areas are known, then the area of the whole is equal to the sum of the areas of the individual parts. It also assumes that areas remain the same when a shape is moved or rotated. If you accept all these assumptions, then you are compelled to accept Pythagoras's Theorem.

This is the essence of mathematics: one lists assumptions (which mathematicians call axioms) and then logically deduces consequences (which are theorems). The less assumptions one uses, and the more unexpected and interesting the consequences one deduces, the more impressive the mathematics. This axiomatic approach was pioneered by the Ancient Greek mathematician Euclid who demonstrated that an enormously rich theory of geometry could be built using just ten axioms. Euclid's axioms are far simpler than the assumptions we made above about area, the first axiom is "Things which are equal to the same thing are also equal to one another". Euclid divided his assumptions into "axioms" and "postulates", but we will call both sets of assumptions "axioms".

Euclid believed his axioms were self-evident truths about the properties of straight lines and angles. However, Euclid's lines are idealised and non-physical and this undermines their claim to being "self-evident". One problem that nagged at mathematicians for centuries was that Euclid assumed the existence of parallel lines, that is lines which can be extended infinitely far without ever meeting. This cannot be tested physically, so how can it be self-evident? For centuries mathematicians tried without success to prove Euclid's assumptions about parallel lines using only the simpler assumptions until it was finally proved that no such proof can exist.

The breakthrough was made independently in the nineteenth century by three mathematicians, Lobachevsky, Bolyai and Gauss. The key observation is to notice that if we replace the word "straight line" in all of Euclid's axioms with "geodesic" (which means the shortest curve on a surface joining two points) then, with the exception of the assumption on parallel lines and the assumption that any line can be extended indefinitely, Euclid's axioms still hold. Since two geodesics on a sphere will always meet, it follows that it must be impossible to show the existence of parallel lines from just the axioms that hold for geodesics on a sphere. It is also possible to construct "negatively curved" geometries similar to that of geodesics on a sphere where all Euclid's axioms hold except the "parallel postulate". Hence that cannot be deduced from the other axioms.

The geometry of geodesics on a sphere is one example of a non-Euclidean geometry. Since spheres are curved, their geometry is different from the geometry of an infinite flat plane. If one studies only small shapes, the results of non-Euclidean geometry only differ minutely from Euclidean results, but as shapes become larger discrepancies begin to appear. We must now be careful to define π as the ratio of the circumference of a Euclidean circle to its diameter. In non-Euclidean geometries this ratio may be slightly smaller or slightly larger than π depending upon whether the surface is positively curved (like a sphere) or negatively curved (like a saddle or Pringle crisp, see Figure 4).

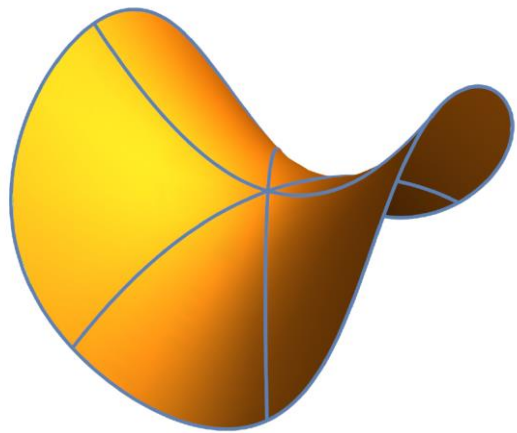

**Figure 4:** A negatively curved "circle". *If a surface is concave or convex in all directions it is said to be positively curved. If it is concave in one direction and convex in another, we say it is said to be negatively curved. All the spokes shown have the same length, r, so we can say the boundary is forms a "circle", but its length is greater than 2πr.*

Non-Euclidean geometry was initially developed out of mathematical curiosity. However, it later formed the basis of Einstein's theory of general relativity, in which space and time are linked together to form a curved space-time. Einstein's theory predicts measurable phenomena such as the bending of light near the sun and the existence of black holes and gravitational waves (Nobel Prize Outreach, 2017). It is no longer possible to reasonably claim that Euclid's axioms are self-evident truths about the geometry of "ideal" space. As a result, many mathematicians today take a *formalist* view of mathematics. In this view, Mādhava's formula is a true statement about the formal theory determined by Euclid's axioms. We don't need to assume it has any connection to physical reality at all.

We can call a collection of axioms and assumptions a mathematical model. Once you have chosen a model, mathematical truths about the model are universal, but the choice of model is determined by the researcher.

How should one choose a model? From the point of view of applied mathematics there are two important criteria. The easier a model is to work with the better, and the more closely it reflects reality the better. These two properties are often in conflict and so some judgement is needed to choose the best model for the situation. Calculations in Newtonian gravity are much easier than in General Relativity and so scientists will often use Newtonian models even though they know they aren't perfectly accurate. As another example, it is useful to model the earth as flat (on a map), even though it is in fact a sphere. A quote often attributed to the Statistician George Box sums up the conventional view of applied mathematical modelling: "All models are wrong, but some are useful."

Physical reality plays a partial role in determining the models used by applied mathematicians, but it is not the full story. The choice of model is also often determined by a desire for simplicity and limitations on the number of meaningfully different but simple mathematical models

Consider as an example the whole numbers. There are many different possible ways of representing a number: one can write numbers in base 10, base 60 or Roman numerals; one can use Chinese or western digits; one can avoid writing and represent numbers as knots on a string, beads on an abacus, lines in sand or fingers and toes. Yet somehow these concepts are all "essentially the same".

To see why let us define a "numbering system" to be any non-repeating sequencing of states or symbols. By a sequence, I mean that it has a first element, an operation "next" acting on elements of the sequence and all elements of the sequence can be found by repeated applications of "next". For example 1,2,3,4,… is a sequence with 1 being the first element and next(1)=2, next(2)=3 and so forth, but i, ii, iii, iv,… is an equally valid sequence.

Suppose you are given two numbering systems, you can match up their elements as follows. First match together their first elements. Next, if x is matched with y, then we will match next(x) with next(y). This completely determines how to map from one numbering system to another. It shows that two numbering systems are essentially the same: the specific labels we use for numbers are completely unimportant.[1]

When two mathematical models are "essentially the same" once one ignores irrelevant details such as notation, we say that the models are "isomorphic". One important activity in mathematics is to classify models up to isomorphism so that we can understand the full range of mathematical possibilities. We have just given an example: up to isomorphism there is only one numbering system. Another somewhat less trivial classification theorem shows that, up to a similar notion of isomorphism, there is only one theory of geometry that satisfies Euclid's axioms. The specific details of Euclid's axioms are not important to determining this geometry, any other set of axioms which identified the same set of geometries up to isomorphism could be used instead. The details used to describe a mathematical model may be somewhat arbitrary and culturally specific, yet the mathematical content of the theory may still be independent of these arbitrary decisions.

---

[1] Notice that this definition of a numbering system has nothing to do with counting! This is an example of the formalist approach to mathematics. We don't try to explain the metaphysical concept of quantity, we simply describe the mechanical aspects ("first" and "next") of a numbering system. However, one can use quantity to create a numbering system: we could use heaps of apples to represent numbers, with "next" corresponding to adding an apple to a heap and "first" corresponding to the "heap" containing only one apple. Since all numbering systems are equivalent, we can use any numbering system to label the heaps of apples. In other words, we can use a numbering system to count apples.

What classification theorems show, is that if one restricts attention to simple, tractable models, then there may be only a limited range of modelling possibilities available. As a result, the same mathematical model can be applied in very different situations. We use the same model of real numbers to model: lengths, volumes, areas, counts, positions in space, durations of time, angles of rotation, bank balances and IQs. In his own work, the first author has examined how the mathematical model Euclid used to describe space can also be used to understand the risk of financial investments (Armstrong, 2018). This doesn't mean that there is some deep connection between the fabric of space and financial risk management, it simply means that the possibilities for tractable mathematical modelling are limited and hence the same model will occur in different places.

In summary, mathematics is used and researched across the globe and has been for millennia. There is strong evidence of some form of "universality" of mathematics as the same highly non-obvious results have been developed independently internationally and can be used to model distinct phenomena.

# Part 2: Decolonising mathematics

Now that we have described the conventional understanding of mathematics, we shall see how this understanding is challenged by the programme of decolonisation. We shall discuss a number of recurring themes, examining both the key texts, and other texts chosen systematically from a sample of the academic literature.

To generate this sample, on 31 December 2022, we conducted a search on Web of Science using the search term "decoloni* AND (math or mathematics or maths)" for matches in any field. This identified 87 texts, which we then filtered to identify those which were written in the English language (eliminating 3 texts) and which discussed the decolonisation of either mathematics, or STEM more broadly (47 texts, mostly on bacterial decolonisation or the decolonisation of other specific STEM disciplines). This left a total of 37 texts which we classified as follows: 8 theoretical papers on decoloniality, 4 papers on critical pedagogy, 12 on teaching indigenous minorities, 7 on teaching in Africa, 2 with particular relevance to the UK, 3 on teaching mathematics to non-native English speakers and 1 on building an online network of researchers.

We use this systematic approach to assess how the decolonisation of mathematics has been understood in the literature.

## Decoloniality

### Key Texts

The theory of decoloniality is a postmodernist critique of the "European paradigm of rational knowledge" (Quijano, 2007, p. 172). Quijano is a key figure in the decolonial programme and states:

> "The colonizers also imposed a mystified image of their own patterns of producing knowledge and meaning. At first, they placed these patterns far out of reach of the dominated. Later, they taught them in a partial and selective way, in order to co-opt some of the dominated into their own power institutions. Then European culture was made seductive: it gave access to power. … Cultural Europeanisation was transformed into an aspiration." (Quijano, 2007, p. 169)

In mathematics, D'Ambrosio (1985) introduced the concept of *ethnomathematics*. He states that "We should not forget that colonialism grew together in a symbiotic relationship with modern science, in particular with mathematics" (d'Ambrosio, 1985, p. 47) D'Ambrosio argues that "[belief] in the universality of mathematics" is becoming hard to sustain in the face of anthropological research and that processes such as counting and ordering are "known" to severally culturally differentiated groups "in a way which is substantially different from the Western or academic way of knowing them" (d'Ambrosio, 1985, p. 46).

> "We know from other sources that, for example, the concept of "the number one" is a quite different concept in the Nyaya-Vaisesika epistemology: 'the number one is eternal in eternal substances, whereas two etc., are always non-eternal'" (d'Ambrosio, 1985, p. 44)

D'Ambrosio argues that "Pure Mathematics" has "obvious political and philosophical ideological overtones" which, for the Third World are "highly artificial and ideologically dangerous" (d'Ambrosio, 1985, p. 47).

In his paper "Western mathematics: the secret weapon of cultural imperialism", Bishop (1990) (who is cited by 12 authors in our sample of 37) argues that "to decontextualise, in order to be able to generalise is at the heart of western mathematics and science; but if your culture encourages you to believe, instead, that everything belongs and exists in its relationship with everything else, then removing it from its context makes it literally meaningless" (Bishop, 1990, p. 57) Bishop states that western mathematics can be distinguished from other forms of mathematics by its values of "rationality", abstraction, "power and control" and "progress" (Bishop, 1990, pp. 56-58). He argues that by promoting these four values, mathematics promotes western "cultural hegemony" (Bishop, 1990, p. 63): "[Mathematics] was (and still is) so clearly useful knowledge, powerful knowledge, and it seduced the majority of peoples who came into contact with it" (Bishop, 1990, p. 58).

Bishop acknowledges that the term "western mathematics" is "in a sense… inapprorpriate", due to the historical and contemporary development of mainstream mathematics outside the West, but then argues that the term is "thoroughly appropriate" since it was "western culture which played such a powerful role in achieving the goals of imperialism" (Bishop, 1990, p. 53). Aikenhead (2017, p. 93) argues that:

> "The Eurocentric impulse to appropriate from other cultures can account for how European mathematicians throughout the centuries seem to have imported ways of mathematizing from earlier cultures but then reconstructed those ideas to fit the European mathematical philosophy or ideology of the time."

### Literature in our sample

Eight of the papers from our sample examine the decolonisation of mathematics from the theoretical lens of the academic theory of decoloniality. Rosa and Mensah (2021, p. 504) argue that "colonialism promoted an epistemicide" and that STEM subjects are a field ruled by "male Eurocentric rationality." Le Roux and Swanson (2021, p. 326) argue that "Scientific knowledge – underpinned by the universal descriptive, analytic and predictive possibilities of measurement, geometry, and methods of probability provided by mathematics – gave legitimacy to [colonialist] hierarchical relations" and this was achieved through "control of subjectivity". Fernandes, Giraldo, and Matos (2022, p. 8) see the goal of decolonising mathematics to be "unveiling how both Mathematics and the patterns of world power cofunction, and creating mechanisms to overcome the myth of Western Modernity". Four of the papers (Rodriguez, 2020, 2021a, 2021b, 2022) defy comprehension: "From the rhizomatic deconstruction as transmethod, the complex objective is fulfilled: to analyze the development of the technologies of the self in Transcomplex Decolonial Mathematics Education."

All of these authors cite Quijano, demonstrating his status as a leading figure in this field.

## Critical Mathematics Pedagogy

### Key Texts

Critical pedagogy is a philosophy of education which has its roots in Freire's "Pedagogy of the Oppressed" (Freire, 1970). Freire proposes that all education is inherently political, and that:

> "The educator has the duty of not being neutral." (Horton, 1990, p. 180)

This is because:

> "Washing one's hands of the conflict between the powerful and the powerless means to side with the powerful, not to be neutral." (Freire, 1985, p. 122)

In critical pedagogy, the purpose of education is to achieve social justice. Freire argues that to achieve this, students must develop "critical consciousness" (conscientização) which will then lead to "the conviction of the oppressed that they must fight for their liberation" (Freire, 1970, p. 67). Freire contrasts this with "the oppressor consciousness" which "tends to transform everything surrounding it into an object of its domination" and which seeks to "in-animate" everything, a tendency he equates with "sadism" and which he views as "necrophilic" (Freire, 1970, p. 59).

> "The investigator who, in the name of scientific objectivity, transforms the organic into something inorganic, what is becoming into what is, life into death, is a person who fears change. He or she sees in change (which is not denied, but neither is it desired) not a sign of life, but a sign of death and decay. He or she does want to study change—but in order to stop it, not in order to stimulate or deepen it. However, in seeing change as a sign of death and in making people the passive objects of investigation in order to arrive at rigid models, one betrays their own character as a killer of life." (Freire, 1970, p. 108)

Freire argues that educators should not follow a "banking" model where students are taught a syllabus determined by the educator. Instead he proposes a "problem-posing" approach where educators help students identify the social justice issues they face and work with them to solve these problems.

Frankenstein (1983) introduced the term "critical mathematics education", linking it directly to Freire's pedagogy. She states:

> "Freire's epistemology is in direct opposition to the positivist paradigm currently dominant in educational theory. Positivists view knowledge as neutral, value-free, and objective, existing totally outside of human consciousness." (Frankenstein, 1983, p. 316)

She emphasizes the potential for statistics to be misleading, and suggests that mathematical skills should be "learned in the context of working on applications which challenge the contradictions involved in supporting hegemonic ideologies" (Frankenstein, 1983, p. 330). For example, she suggests teaching division by getting students to show that "each of the richest 160,000 taxpayers got nine times as much money as the maximum AFDC grant for a family of four" (Frankenstein, 1983, p. 327). This is necessary to avoid supporting hegemony as in conventional mathematics education

"Even trivial math applications like totaling grocery bills carry the ideological message that paying for food is natural" (Frankenstein, 1983, p. 328).

Skovsmose (1994) provides a specific case-study of critical mathematics pedagogy. He argues that "mathematics is formatting society" (Skovsmose, 1994, p. 36), and argues, citing (Davis & Hersh, 2005, pp. 120-121), that this can be seen through its use in rules such as "income tax", "taking a number in a bake shop", and similar systems where mathematical modelling is used to prescribe unambiguous rules "often for reasons known only to a few" . To educate students in this, he conducted a series of lessons where student were asked to design their own child benefit rules. He discusses the outcomes of this project in terms of the extent to which students were able to perceive the way mathematics formatted society, but not in terms of gains in traditional mathematical knowledge by the students.

Public understanding of critical mathematics pedagogy has come, not through the works of these authors directly but instead through media reporting such as (Gearty, 2017; Richardson, 2017). This article discussed the paper (Gutiérrez, 2017a), highlighting quotations such as "On many levels, mathematics itself operates as Whiteness" (Gutiérrez, 2017a, p. 17) and her concern about minority students who "have experienced microaggressions from participating in math classrooms… [where people are] judged by whether they can reason abstractly" (Gutiérrez, 2017a, p. 18). In Gutiérrez (2017b) defends her approach and describes its academic origins in the work of scholars such as Bishop, Frankenstein and Skovsmose.

## Literature in our sample

Four texts in our sample focussed on critical pedagogy. Bhattacharya et al. (2022) report a "decolonial" approach to creating a "College Math Academy" teaching precalculus mathematics using project-based learning in a manner that would provide opportunities to "critique dominant power structures" (Bhattacharya et al., 2022, p. 235). The students quoted reported increased awareness of social justice issues, but do not mention any specific mathematical skills they developed. De Roock and Baildon (2019, p. 285) examine an online portal designed to promote STEM education in Singapore and find that it "leverages student interests, self-discovery and diverse identities to recruit them into (neoliberal) figured words of learning". Nicol (2016) gives an example of how students were able to use mathematics to compare indigenous logging practices on Haida Gwaii (a Canadian archipelago) with European practices. Nicol argues that this experience gives a counter narrative to claims that critical pedagogy can lead to "brainwashing" "as neither teachers nor students came to think in one homogenous perspective" (Nicol, 2016, p. 437), though she provides no details of the perspectives held before or after the intervention, or the impact upon students' mathematical understanding. Trinder and Larnell (2015, p. 220) outline a proposal for a new theoretical framework for mathematical pedagogy which they hope will "reinvigorate the central question of how knowledge(s) for teaching mathematics will be regarded, whose knowledge(s) will be included, and how we can do this work in ways that subvert the imperialistic tendencies of education."

## Indigenous Mathematics Education

### Key Texts

In countries such as Canada, New Zealand and Australia with minority indigenous populations, the concept of "decolonisation" could be understood as ensuring equitable educational outcomes for this indigenous population.

The Truth and Reconciliation Commission of Canada (2012) Calls to Action identify seven points pertaining to education, the tenth of which includes improving attainment levels for indigenous children and providing culturally appropriate curricula (though as of June 2022 no new legislation was being drafted (CBC News, 2022)). This seems to have little relation to the postmodern theories of decoloniality, unless one changes the metric of attainment to shift away from traditional metrics or one defines cultural appropriateness to require fundamental epistemic change. Nevertheless, the advocates of indigenizing the science curriculum found in (Widdowson, 2016) do appear to take inspiration from postmodernism.

The New Zealand National Curriculum is explicit in requiring "Equal status for mātauranga Māori [Maori knowledge]" (NCEA, 2021), which accords well with the postmodern theory. The Biology and Chemsistry syllabus required teaching the atomic theory using concepts such as "mauri" or life force (Matzke, 2022), though these materials have since been changed after protests from scientists.

## Literature in our sample

Eleven out of the twelve texts in this category sought to question the epistemic privilege of modern science, frequently characterising this as Western. Nicol, Gerofsky, Nolan, Francis, and Fritzlan (2020, p. 192) observes "Nearly Universal And Conventional mathematics offers a particular worldview that values rationalism, objectivism, power, control and progress". Garcia-Olp, Nelson, and Saiz (2022, p. 4) propose a programme which is "not just centering Indigenous knowledge" but which "live[s] out Indigeneity" through the incorporation of concepts such as "Heart Work" and "The Tree of Life" (Garcia-Olp et al., 2022, p. 2). Baker et al. (2021, p. 91) advocates "engagement that results in families seeing their epistemology and values privileged". Kulago, Wapeemukwa, Guernsey, and Black (2021, p. 346) state: "the purpose of decolonial curricula is to rupture the epistemic barriers of Western knowledge". Kearns, Tompkins, and Lunney Borden (2018, p. 236) "use indigenous ways of knowing, being and doing to challenge mainstream understandings of teaching and learning." Chinn (2007, p. 1250) asserts: "Western science methods of knowledge building … are antithetical to a Hawaiian world view". Velez et al. (2022, p. 1) questions what it might mean to improve outcomes for indigenous students and states that "a paradigm shift from Western evaluation methodologies to Indigenous evaluation is necessary when evaluating STEM programs".  Harper (2017, p. 359) discusses a "cross-cultural learning community that privileges Karen cultural knowledge". R. Eglash et al. (2020) measure the success of their intervention in terms of the extent to which students develop an appreciation of Anishinaabe ethnomathematics. After the intervention one student commented: "Going to Mars is a concept that scientists have aimed for years and could be made possible by the Anishinaabe traditions" (R. Eglash et al., 2020, p. 1584). The same intervention is discussed in Ron Eglash, Bennett, Drazan, Lachney, and Babbitt (2017).

The only text on educating indigenous minorities that was unclear in its epistemological position was Sarra and Ewing (2014). They cite Bishop when discussing the importance of avoiding "enculturation", but the intervention they describe appears to consist of quite conventional pre-school classes aimed at aboriginal children.

Culturally appropriate education was also a common theme. Garcia-Olp et al. (2022) suggest having students compute the volume of tipi. Ron Eglash et al. (2017, p. 1573) describe a number of approaches ranging from "a unit called 'wigwametry' using circular model wiigiwaam construction to allow Native students to investigate the value of the mathematical constant Pi" to lessons drawing analogies between Anishinaabe boatbuilding and Bezier curves. Aikenhead (2017, p. 81) discusses the risk of tokenism and stereotyping in decolonisation as follows:

'Sterenberg (2013a) repeated an insulting word problem "Imagine a band of 250 Aboriginal People. Each tipi can hold approximately eight people. Calculate how many tipis would be needed to house the entire band" (p. 21). This is insulting because Indigenous people would not divide themselves in the hypothetical way stated in the word problem. Relational and spiritual factors would dominate. And the required hypothetical state of mind itself, captured in the word imagine, is a value embraced strongly in the culture of school mathematics but would generally be foreign to an Indigenous culture in the context of people choosing a tipi to enter, because it conflicts with how to live in a good way.'

Not all of the authors of Nicol et al. (2020, p. 201) were in full agreement. Francis reports successfully using robotics to inspire her indigenous students and states: "robots are NOT place based or culturally responsive".

Critical pedagogy was also a recurring theme. However, interventions were typically co-created in conjunction with community elders as well as students, e.g. (Harper, 2017).

## African Mathematics Education

### Key Texts

In his review of post-colonial teaching of STEM in Africa, Yoloye (1998) reports on a number of major international conferences that emphasized that the development of STEM capacity should be a major priority in developing nations. This should include "reform[ing] the content of education in the areas of curriculum, text books and methods, so as to take account of the African environment" (UNESCO, 1961). Yoloye describes a number of interventions that sought to provide a science education to poor communities within Africa, many involving collaboration with Western scientists. Yoloye does not mention "decolonisation" or any of the thinkers discussed above, and these conferences predate the development of the theory of decoloniality, so we might term these approaches "mainstream science education". Although (Crowell, 2023) considers the African Institute of Mathematical Sciences to be an outgrowth of the movement to decolonise mathematics, it was founded in 2003 and makes heavy use of collaborations with Western mathematicians so it seems more likely to us that it grew out of the mainstream tradition, a view supported by (Hutton, 2023).

It is, however, possible to detect a decolonial theme in science policy in South Africa as it emphasises the importance of valuing "indigenous knowledge systems" (Department of Science and Technology, 2005). Mudaly (2018) (in our sample) describes this as decolonisation, but is critical of the lack of specifics in the policy which he argues makes it difficult to implement.

### Literature in our sample

A common theme in the African literature is that of determining the boundaries of the decolonisation programme and four out of seven authors raised this. Singh-Pillay, Alant, and Nwokocha (2017, pp. 131-132) examine attitudes of science teachers in Nigeria to the concept of "indigenous knowledge" and finds that "teachers conceive [indigenous knowledge] as technological knowledge with a science base" and "do not embrace the metaphysical or spiritual aspects". Taylor and Cameron (2016, p. 42) recommend teaching that considers the intersection of indigenous knowledge with scientific knowledge, thereby acknowledging that not all indigenous knowledge is "suitable for use in the classroom". De Beer (2016, pp. 48-49) advocates for the inclusion of indigenous knowledge in the curriculum because it may lead to "hopefully an improvement in cognitive skills and better test scores", but he notes that he does not wish to propagate "pseudo-scientific approaches" or approach the topic from the political perspective of "decolonising the curriculum". Schubring (2021, p. 1465) studies the history of STEM education, finding clear examples of culturally inappropriate material in historic textbooks. However, Schubring does question some aspects of the theory of decoloniality, remaking that "Walsh's (Mignolo & Walsh, 2018) programmatic claim that decoloniality should 'displace Western rationality as the only framework and possibility of existence, analysis and thought' is quite strong" (Schubring, 2021, p. 1465).

Mudaly (2018, p. 72) argues that South African policy needs to provide more examples of how to implement decolonisation in practice. He found that "practicing teachers, including some from African and Indian communities, could not identify aspects of a colonised curriculum, so they did not know how they could easily decolonise it. Many of their responses were superficial and they provided examples and solutions that appeared to be coerced and contrived… Whenever they attempted to illustrate the use of indigenous knowledge in mathematics, the only example they chose was the cylindrical hut with a conical thatched roof."(Mudaly, 2018, pp. 72-73)

While the papers above suggest a common understanding of decolonisation in the African literature there were two outlying papers. Khoza and Biyela (2020, p. 2666) define decolonisation of the curriculum very broadly as "critiquing and renewing the knowledge of curriculum" and we could not see any relationship between the intervention they describe and colonisation. We were unable to understand Liccardo (2018, p. 12) who "proposes an infinity symbol as a theoretical-empirical framework to illustrate the psychosocial rhythms of inclusion-exclusion and the reproductions of disciplinary and institutional cultures at historically white universities".

## UK Mathematics Education

### Key Texts

The programme of decolonising mathematics in the UK was discussed by Brodie (2016) on UK academic forum The Conversation. She argues that "universities must do explicit identity work with their students" and highlights the role of ethnomathematics and critical mathematics. In 2022 a draft document by the Quality Assurance Association for Higher Education (QAA) recommended that all UK mathematics degrees teach a "decolonised view" of mathematics (White, 2022). The final document (QAA, 2023) contained different wording but still associated decolonisation with teaching about "problematic" aspects of the history of mathematics such as the existence of Nazi mathematicians and statisticians who have supported eugenics.

At Durham University lecturers were asked to consider giving short biographies of the mathematicians whose work they present and if they are "almost entirely (or even completely) white and/or male, ask yourself why this" (Sommerville, 2022). Speaking at a London Mathematical Society education day, Crossman and Ogundimu (2023), both from Durham, defined decolonisation as:

> "[A]n ongoing process of decoloniality, which deliberately diminishes predominant voices, disinvests from power structures, devalues hierarchies, decentres knowledge production, and diversifies ways of knowing. Its movement is to generate empathy and mutuality towards entangled subjectivities, such that individuals become many-sided, social beings, capable of lifelong transformation" (Hall, 2020)

Crossman did emphasize in his talk that "1+1=2", but Brigitte Stenhouse, speaking later stated that she did not agree with this, showing a slide of a statement by the University of Durham saying that "2+2=4" and arguing that in her view we should also consider other ways of knowing. Stenhouse appears to adopt a different position in (Barrow-Green & Stenhouse, 2022).

Both pairs of speakers broadly advocated increasing the emphasis on non-white and non-European contributions to mathematics. A project at Queen Mary University of London (QMUL) produced a list of biographies of mathematicians from minority backgrounds that was intended to be used as a resource in teaching (Onus et al., 2022) and which was promoted at King's College London as a resource on decolonising mathematics. It lists mathematicians such as Tyler Kelly, a non-binary algebraic geometer at the University of Birmingham and Alex Fink, a queer combinatorics researcher at QMUL alongside better known names such as the Islamic mathematician Al-Khwarizmi from whose name the word algorithm derives. This resource lists four Jewish mathematicians all of whom have some other characteristic (female, disability, "queer") and who appear to have been selected on this basis. As a result, figures such as Einstein are omitted. Another notable omission is Leonard Euler who was described by Laplace as "the master of us all" and who in 1775 wrote on average one mathematical paper per week despite total blindness (Finkel, 1897).

Borovik (2023) also finds that the decolonisation of mathematics in the UK focusses upon teaching the history of mathematics. Borovik is critical of the decolonial programme and illustrates how the history of mathematics might be taught in a manner that contradicts the theory of decoloniality, using the example of linear algebra.

## Literature in our sample

C. K. Raju (2017) has made two prominent contributions to UK understanding of the decolonisation of mathematics, though his views are not specifically about the UK.  Firstly, (C K  Raju, 2016) was featured on the UK academic forum The Conversation but was subsequently deleted after questions were raised about the quality of Raju's scholarship, something that he and others have subsequently attributed to racism. Secondly, Raju's approach to decolonising mathematics provides one of the most concrete proposals in the collection "Rhodes Must Fall" (Kwoba, Chantiluke, & Nkopo, 2018). While Raju's thinking has some similarities with scholars such as Quijano, his philosophy of mathematics and mathematics education is distinct. Raju's work has a broad scope, but as one example he claims to be able to teach students calculus in five days (C. K. Raju, 2017). This is a remarkable claim as traditional calculus courses may last a year or more. However, one looks into the detail of what Raju is teaching (C K Raju, 2013), he is only teaching students to use a computer program that is able to find approximate solutions to a certain class of calculus problems numerically, namely a numerical differential-equation solver. He argues that the "calculation of symbolically complex derivatives and integrals" which is taught in traditional courses, "is completely pointless today, when it can be done in a jiffy using open-source symbolic manipulation programs such as MAXIMA" (C K Raju, 2013, pp. 4-5). Raju's claim is therefore broadly analogous to saying that one can teach foreign languages in just a few hours by showing students how to use automatic translation software.

The other work in our sample that was of particular relevance to the UK context is (Skopec et al., 2021). This article comments "there are few examples of decolonization efforts in … STEMM disciplines in the UK and no empirical studies to our knowledge" (Skopec et al., 2021, p. 2), though they sought to address this using a mixed-methods methodology including an "implicit association test" looking at attitudes to research from lower and middle income countries. They do not provide sufficient details of their test for me to understand what it entailed. Their principal theoretical contribution is that resistance to decolonisation within STEMM may arise from "epistemic fragility", which they relate to the concept of "white fragility" introduced by DiAngelo (2011).

## Miscellaneous Literature in our Sample

Bennett (2016) describes the use of computer tools that combine art with STEM disciplines with a focus on black identity, for example a program that allowed students to develop designs based on cornrow hairstyles. They report statistically significant improvements in mathematics outcomes compare to a control group, but do not compare their intervention with similar interventions without the focus on identity.

Three texts examined the challenges of teaching mathematics to non-native English speakers (Dlamini & Kamwendo, 2018; Parra & Trinick, 2018; Ruef, Jacob, Walker, & Beavert, 2020).

One text developed an online community which "drew on speculative fiction to define a #VanguardSTEM hyperspace as a fluid "place-time" that is born digital and enabled by social media, but materializes in the physical world for specific purposes" (Isler et al., 2021, p. 1).

## What does it mean to decolonise mathematics?

Having examined the literature, we are now able to answer the question of what it means to decolonise mathematics according to its advocates. The recurring theme is that we should value other "ways of knowing" equally, or perhaps even more highly, than scientific and mathematical approaches to knowledge.

What constitutes a valid "way of knowing" is not discussed in the literature in our sample. This perhaps reflects the idea that one should not attempt to constrain what is and what is not legitimate knowledge. In the introduction to the collection "Ways of Knowing" Harris (2007, p. 20) traces the origin of the concept of different ways of knowing to Giambattista Vico (1668-1744) explaining: "Reason, for Vico, was only one way to comprehend the human condition and discover truth: there were others such as imagination and empathy (called 'kinds of knowing' earlier)". Vico pioneered the social constructivist approach to science that argues that scientific knowledge is a social construction rather than an objective analysis of the external world, stating "*Verum esse ipsum factum*" ("truth is itself something made"). In their book "The Social Construction of Reality", Berger and Luckmann (1991, p. 83) propose the following broad concept of knowledge: "theoretical knowledge is only a small and by no means the most important part of what passes for knowledge in a society… It is the sum total of 'what everybody knows' about a social world, an assemblage of maxims, morals, proverbial nuggets of wisdom, values and beliefs, myths, and so forth…" A similar concept of knowledge can now be found in school syllabuses: (Popov, 2020) lists various "ways of knowing" including sense perception, emotion, imagination, intuition, faith, memory and reason.

This understanding of different "ways of knowing" appears to have become standard in anthropology. However, within the sciences the term "knowledge" has a conflicting meaning as "justified true belief". We believe this concept remains a useful one and so for the remainder of this paper we will use the word "knowledge" solely to refer to justified true beliefs and "beliefs" for all other beliefs.

This allows us to state the conclusion of our literature review more clearly. The central recurring theme in the English-language literature on decolonising mathematics is to question the epistemic privilege of scientific and mathematical reasoning over other beliefs. This accords with the well-studied theory of decoloniality, and more generally with postmodernist epistemology. A consensus emerges from the literature that this is the most generally accepted meaning of decolonising mathematics.

Many of the authors we have studied would go much further. Some would want mathematics to be taught in a way that promotes an anti-capitalist account of social justice. Others would want to teach that mathematics was created for the purpose of European domination. We cannot argue that these ideas are essential to the decolonisation programme as not all the authors we studied expressed these views directly. However, it is interesting that the non-African authors in our sample did not choose to critique such views either.

In the African literature a number of authors argue that mathematics education should only consider indigenous beliefs which are justified and true. These authors support incorporating such traditional indigenous knowledge into the mathematics curriculum. For these authors, this is what distinguishes decolonised mathematics from other approaches to culturally appropriate mathematics education. These authors also discuss the attitudes of other African mathematicians, many of whom it seems have some doubts about the decolonisation programme, and  The journal Nature has been criticised for an apparent publication bias in favour of "critical social justice" approaches (Abbot et al., 2023), which include decolonisation. Our study suggests that such a bias would discriminate against African authors. Táíwò (2022) critiques decolonisation for not "recognizing or engaging with African intellectual discourses" that do not conform to a decolonial ideology.

Within the UK, decolonisation of mathematics has taken a specific form that is quite different from that of the international literature. Proponents of decolonising mathematics in the UK do not propose that students should study the ethnomathematics of stone circles or early Christian astronomy. Nor is anyone seeking to inspire students by having them compute the circumference of King Arthur's Round Table. The UK programme of decolonisation of mathematics questions the epistemic privilege of mathematical reasoning by proposing that much of mathematics was developed for racist ends and that the contemporary mathematics curriculum excludes contributions from minority groups due to ingrained historical prejudice.

Our literature review bears out the view expressed in Nature (2023) that epistemology is a central theme of decolonisation. But this is not the end of the story: having identified that decolonisation rejects scientific epistemology, it appears that the rationalist approach we have been taking may have been fundamentally misguided. Perhaps it would be better to try to understand decolonisation emotionally. Indeed, perhaps it is foolish to seek a purely rational explanation of why the murder of George Floyd by a policeman in Minnesota should lead to changes in the mathematics curriculum at the University of Oxford (Woolcock, 2020). Emotionally however, the connection is understandable: white people in the UK have a collective guilt for the slave trade and so they must all do something, whatever it may be, to assuage this collective guilt, and this includes the mathematicians.

While epistemology is the recurring intellectual motif of the literature we reviewed, the theme of decolonisation as a moral imperative is the recurring emotional theme. Biggs (2020) remarks that the book Rhodes Must Fall (Kwoba et al., 2018) is "surprisingly vague on specifics", yet this does not dampen the authors' certainty that decolonisation, whatever form it may need to take, is essential. Haxton, Karodia, Taylor, Williams, and Lalemi (2022) explain that "[decolonisation of the curriculum] remains poorly understood within the academic science community" but nevertheless found "Staff were particularly keen to find out what they could do to decolonise their curricula". The emotional, spiritual and moral necessity of decolonisation are highlighted by a number of authors in our sample. Kearns et al. (2018, p. 249) state "We are reminded that decolonizing work has a spiritual dimension to it, as Mi'kmaw Chief PJ Prosper so eloquently stated, 'moving from the head to the heart'". Aikenhead (2017, p. 78) observes "I have come to realize that for my contribution to go above tokenistic culturally relevant content, it is not only the content that must be decolonised but also myself". Fernandes et al. (2022, p. 7) say they "understand decoloniality as a choice: a decision to share a political agenda engaged in struggle, resistance and insurgency against the various traces and effects of coloniality we are traversed by".

Linking the intellectual and emotional theme is the idea (sometimes implicit, sometimes explicit) that privileging mathematical knowledge over other beliefs is racist. This explains the central importance of the social construction of knowledge to decoloniality: if one accepts a universalist view of mathematics, there is nothing racist about privileging mathematics. The key authors we have discussed who explicitly address mathematics are Bishop, d'Ambrosio, Frankenstein and Skovsmose. These authors all succeed in showing that some of the superficial aspects of mathematics (such as choices of units) are social constructed, as are many of its applications (such as systems of taxation). However, none of this contradicts the universalist view we have put forward. This universalist view acknowledges that mathematics is a human endeavour, and that modelling is subjective, but argues that nevertheless it is still possible to identify universal mathematical truths. Different societies may independently develop identical truths. Bishop is the only author we studied who provides a clear counter to this view, but to do this he must argue that rationality and abstraction are western values. If one is willing to reject rationality and abstraction as universal, we agree it is no longer possible to argue that mathematics is universal. However, without rationality it seems difficult to argue cogently at all. Perhaps all one can do is look into one's heart and ask, which is more racist: the view that mathematical knowledge is more reliable than other beliefs or the view that rationality is a Western colonialist trope?

# Part 3: Is the culture of mathematics racially discriminatory?

One can believe that the mathematics community is, or was, systemically racist whether or not one agrees with decolonialist critiques of the traditional epistemology of mathematics.

In this section we will look for evidence to identify whether there is evidence of systemic racism within mathematics. We focus on the experience in the UK and apply a primarily empirical methodology.

## Demographic Data

It is natural to ask whether ethnic minorities are underrepresented at different levels of the education system. Although the London Mathematical Society has produced reports on the representation of women in mathematics (Ortus Economic Research, 2023), they do not appear to have published similar data on race since 2013 (London Mathematical Society, 2013). Indeed, the LMS women in mathematics committee was only renamed to the Committee for Women and Diversity in Mathematics in 2020.

In Figure 5 we have plotted the proportion of UK-national academic staff and professors in the UK mathematics who are Black, Asian or Minority Ethnic (BAME) by age. We have compared this to the proportion of degree holders (irrespective of subject) in England and Wales who are BAME, obtained using census data. We have only shown academic staff who are UK-nationals to try and minimize the issues caused by the legal requirement before Brexit to give priority to candidates from Europe. However, as staff may have then chosen to adopted UK-nationality, this adjustment to the data will be imperfect, It also makes little visual difference to the resulting chart.

This plot suggests that BAME staff may get promoted later than white colleagues and that it may historically have been more difficult for BAME staff to obtain an academic job. However, there is no indication of bias in the hiring of staff under the age of 40.

Our analysis of this quantitative data is somewhat crude. A detailed quantitative analysis of demographic data in mathematics is surely overdue but is beyond the scope of this paper.

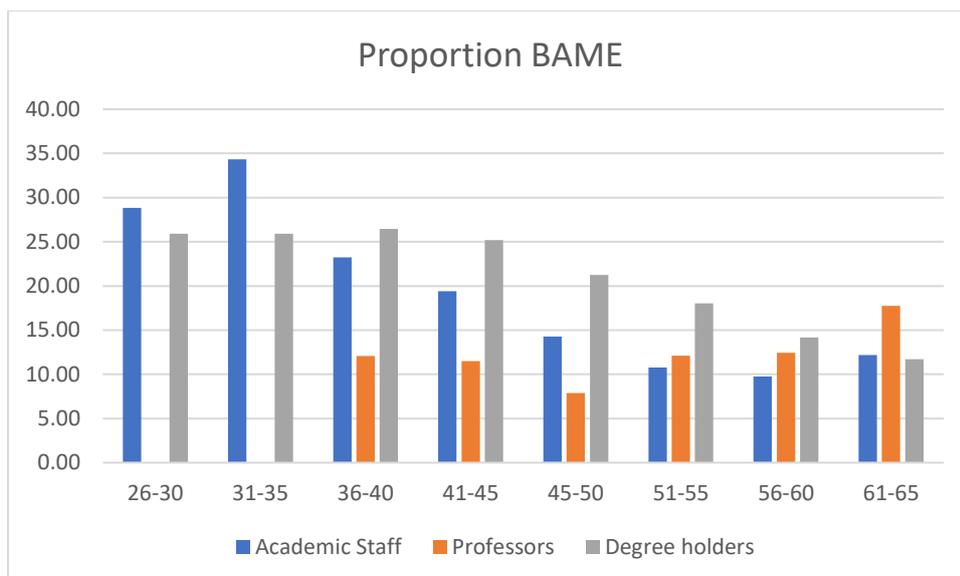

**Figure 5:** Proportions of UK-national professors and academic staff in mathematics who are BAME by age using HESA data for 2021/22 (Higher Education Statistics Agency, 2022). Proportion of degree holders (any subject) who are BAME using census data for 2022 (Office for National Statistics, 2023a), filtered for highest qualification of degree or above. Individuals whose ethnicity is unknown are not included in the data. Data on professors under the age of 36 is supressed in accordance with the HESA rounding methodology. Note that census data is grouped by age categories that start one year earlier: 25-29, 30-34 and so forth.

## The history of mathematics

Many of the authors in our literature review suggested that: mathematics education paints a false history of mathematics as a western creation; mathematics, and especially pure mathematics, is a colonialist creation; or both.

For example, d'Ambrosio (1985, p. 47) states that pure mathematics "is the natural result of the evolution of the discipline within a social, economic and cultural atmosphere which cannot be disengaged from the main expectations of a certain historical moment", by which he means the period of European colonialism. In saying this, d'Ambrosio ignores the Ancient Greek contribution to mathematics. In addition, he is simply asserting that because the development of (some) mathematics and the development of colonialism were contemporaneous they must be linked.

Other authors give more detailed examples. Prescod-Weinstein (2020, p. 427) gives the following examples:

> "Science, mathematics, and slavery were intimately connected: whether it was the early evolution of insurance and actuarial science to calculate the value of jettisoned cargo—brutally murdered people—or efforts to minimize the bow wave—the wake—of ships, to make them faster, to speed the movement of kidnapped Africans from the torturous Middle Passage to a tortured lifetime and usually death in the bondage of chattel slavery"

However, this seems more to establish the undoubted horrors of slavery than the existence of an intimate relationship between mathematics and slavery. The QAA (2023, pp. 6-7) give the following examples, which they recommend should be taught to all university mathematics students:

> "some early ideas in statistics were motivated by their proposers' support for eugenics, some astronomical data were collected on plantations by enslaved people, and, historically, some mathematicians have recorded racist or fascist views or connections to groups such as the Nazis."

As the repetition of the word "some" here suggests, the QAA leave out a great deal. MacLane (1995) describes how the Mathematical Institute in Gottingen was decimated by Nazi antisemitism, aided by student protests. Mathematicians working at Bletchley Park led by Alan Turing played a significant role in the Nazis' defeat. While the QAA emphasizes the aspects of the history of mathematics that support the decolonial thesis, those aspects of mathematics that contradict it go unmentioned, and it seems that in the QAA's view, these do not form an essential part of a mathematician's education.

Barrow-Green and Stenhouse (2020) give another example:

"Mathematics is key to solving the longitude problem in the 18th century, which enabled mariners to locate themselves accurately at sea and gain dominance over the waves, which was essential to Britain and its imperial and colonial aims."

However, it seems one could make an equally convincing argument that ship's biscuits were essential to Britain's colonial aims. Indeed, the sea itself appears equally culpable. Knowledge of mathematics, like any resource, may be used for good or ill.

Having reviewed these examples of mathematics' colonialist credentials, it is worth considering the large amount of mathematics which does not appear to have been motivated by colonialism. Central topics in mathematics such as the theory of gravitation, the calculus, non-Euclidean geometry, number theory and so forth do not have any obvious relationship with colonialism, and it seems advocates of decoloniality have not been able to identify such relationships either.

In addition, some European mathematicians have been notably open to non-Western ways of thinking. Medieval European mathematicians such as Fibonacci worked to popularise mathematics of Indian, Arab and Ancient Greek mathematicians.  As a result, Western mathematics has its origins in a project not entirely dissimilar to ethnomathematics in its openness to other cultures, but diametrically opposed in its preference for rational argument over one's own cultural values.

We have already discussed the decolonialist view that rationality is a Western value, not shared by other cultures, and noted that the international history of mathematics appears to contradict that. We should also note that Western society is by no means uniformly rational. European colonialists had their own spiritual views and sought to spread Christianity. Christianity has a history of opposition to science: Galileo was placed under house arrest by Pope Paul V for advocating for the theory of heliocentrism; Darwin's theory of evolution, was, and still is, rejected by many Christians (Wilberforce, 1860). The role of Christianity in colonialism is largely undisputed, and its influence can still be seen in the residential-school scandal in Canada (Pope Francis, 2022). While many authors on decoloniality consider Ancient Greek mathematics to be Western, this ignores the fact that much of the greatest work in Greek mathematics took place in Alexandria, in present-day Egypt, and that the Ancient Greeks were not Christians.

A more compelling case for racism in the history of mathematics is the naming of mathematical results. As we have seen, Mādhava was the first to discover the Mādhava-Leibniz formula, yet it is mostly called Leibniz's formula by mathematicians. Similarly, the Japanese mathematician Seki Takakazu (c1642-1708) discovered the Bernoulli numbers independently of Bernoulli, and Seki's results were published a year earlier than Bernoulli's (both were published posthumously). Nevertheless, they are typically called Bernoulli numbers.

One example that has gained prominence is the Chinese Remainder Theorem. This theorem has a complex history. It has been studied by a number of mathematicians, but its first known statement is as a special case with specific numbers by Chinese mathematician Sunzi and the first full version with proof by Qin Jiushao. Dickson (1926) gives the first recorded use of the name "Chinese Remainder Theorem". One might see this name as respectfully acknowledging the complex provenance of the theorem, but it has been argued that naming the theorem after a country rather than an individual is racist. Kim (2021) argues this with considerable force:

> "Why did Dickson remove Sun Tzu's name from the theorem? We can't know what was in his heart, but we know that Dickson made the choice amid a surge of anti-Asian violence in the United States stretching back to the late-19th century. For example, in Rock Springs, Wyo., in 1885, a White mob torched the local Chinatown and killed 28 Chinese immigrants. Two years later, in the Snake River Massacre, the mutilated bodies of 34 Chinese miners were found floating down a river in Oregon, butchered by White miners upstream."

What Dickson said in full was that one line in his proof followed "by the Chinese Remainder Theorem". It interesting to note that mathematicians also refer to "Polish spaces" in order to acknowledge the contribution of a number of Polish topologists, but we are not aware of this eliciting similar objections.

Mathematicians frequently misattribute or fail to attribute concepts. Fermat's Last Theorem was famously not proved by Fermat, but by Andrew Wiles centuries later. The Fundamental Theorem of Algebra does not attribute Gauss, and the fundamental theorem of calculus does not attribute Gregory. Seki's work was not known in the West on publication because of Japan's isolation policy of Sakoku.

The work of Mādhava and the Kerala school has taken a long time to be appreciated in the West. Roy (1990) documents how Whish reported these works in a paper of 1835 which had negligible impact, until Rajagopal and associates began to raise the profile of the work. The Kerala school was not mentioned in the 1968 or 1991 editions of the text book (Boyer & Merzbach, 2011), only finally appearing in the third edition of 2011. The first general history of mathematics textbook we found on the open shelves at University College London library that mentioned the Kerala school was (Hodgkin, 2005). Almeida and Joseph (2007) suggest this might be attributed to Eurocentrism with origins in colonialism. It certainly does not appear to reflect well on historians of mathematics that it has taken so long to become known.

Knowledge of the Kerala school has yet to filter into the mathematics community: an informal survey of mathematics academics at King's College London showed that 0 out of 12 academics surveyed knew who Mādhava was. This is very likely to be a consequence of a widespread lack of interest among mathematicians in the subject's history. Only 5 academics in the same survey were aware that Euler was blind in later life despite Euler's central role in the development of calculus. For many mathematicians, what is interesting is why a result is true rather than who proved it. We find this understandable and do not seek to argue that all mathematicians should engage in an extensive study of the subject's history.

## GCSE mathematics textbooks

We can apply the Durham decolonisation recommendations to the UK GCSE mathematics syllabus (for schoolchildren aged 14-16) and identify the mathematician associated with the introduction of each concept. One can debate who best to choose, but our suggestions are below:

1. Number: Aryabhata

2. Algebra: Al-Khwarizmi
3. Ratio, proportions and rates of change: Pythagoras
4. Geometry and Measure: Euclid, Pythagoras
5. Probability: Pascal
6. Statistics: We feel there is no obvious choice at GCSE, perhaps one might choose to discuss John Snow and Florence Nightingale for their medical applications of statistics.

This is an ethnically diverse list, but lacks gender diversity. The selection of topics appears to be based on the fact that the material is both elementary and important, providing an essential foundation in mathematics.

We examined two GCSE mathematics textbooks (Fearnley et al., 2015; Morrison, Smith, McLean, Horsman, & Asker, 2015). We first looked to see whether the names of individuals used in problems were diverse. 11% of given names used in problems were not of North American or European origin, whereas 18% of the population identifies as BAME (Office for National Statistics, 2023b), though many of them have given names of North American or European origin. 30% of the cartoon avatars that appeared in the AQA textbook were BAME.

We also examined the coverage of history of mathematics in GCSE textbooks. History is only discussed fleetingly, with 38 references in total. 15 are references to Christian European history, 11 to Classical Greek and Roman history, 10 to international history or prehistory, and 2 were not sufficiently specific to classify. The strongest acknowledgement of mathematics international history in these textbooks was:

> "Mesopotamian, Chinese and Indian mathematicians all independently discovered Pythagoras's theorem. However, the Greek Pythagoras ended up getting all the credit" (Morrison et al., 2015)

However, there is an asymmetry in which mathematicians are named. Pythagoras, Euclid and Pascal are all named in the textbooks, but Aryabhata and Al-Khwarizmi are not.

The history of mathematics is taught only as an aside within the textbooks and is not examinable.

We conducted three interviews with schoolteachers to gather their views on the teaching of the history of mathematics. They all said that they did not directly teach the history of mathematics, identifying the need to teach core syllabus content as the primary reason. They all agreed it would be potentially useful to teach the history of mathematics, but one teacher questioned whether it would be worth the effort as students might not be particularly interested.

At one of the schools we visited there were posters about mathematicians on the doors of the mathematics teaching rooms. These focussed on mathematicians of diverse identities, but were biased towards showing US mathematicians.

We report fuller details of this study in (TODO cite India's preprimt)

## A Masters Course

In a module taught by the first author at King's College London, each week short biographies are given of the mathematicians behind each theory. Eleven mathematicians are listed. Seven are white men with Jewish backgrounds (Markowitz, Wiener, Einstein, Scholes, Feynman, Kac, Lax). Three are white non-Jewish men (Merton, Brown, Black). One is a Japanese man (Ito).

The absence of women might be attributed to the historical underrepresentation of women in science. The "overrepresentation" of Jewish mathematicians is harder to explain. Gerstl (2014) suggests "Jewish cultural values" such as "dedication to education". Whatever the cause, Jewish mathematicians have historically faced considerable antisemitism, not only in Nazi Germany but also within the US mathematical establishment (Nadis & Yau, 2013).

## Cinema

Wikipedia contains a list of 17 biographical films about mathematicians. After removing the fictional film Enigma, the demographic characteristics of the mathematicians are shown in Table 1. 7 out the remaining 16 films are about white male mathematicians. If one excludes disabled and gay mathematicians and those of Jewish origin this leaves only one film, about Descartes. It is perhaps not surprising that cinema provides a skewed account that focusses on personal identity as its aim is entertainment rather than scholarship.

|  | Number of films | Mathematicians |
|---|---|---|
| Non-white Female | 2 | Sun-Yung Alice Chang, Fan Chung, Katherine Johnson, Dorothy Vaughan, and Mary Jackson |
| White female | 4 | Hypatia, Sofya Kovalevskaya, Maryam Mirzakhani |
| Non-white male | 3 | Srinivasa Ramanujan, Yitang Zhang |
| White male with disability | 3 | Stephen Hawking, John Nash |
| Gay white male | 1 | Alan Turing |
| White male of Jewish origin | 2 | Richard Feynman, Paul Erdős |
| Other white male | 1 | Rene Descartes |

**Table 1:** Demographic characteristics of mathematicians in biopics, using US census definition of white. Source (Wikipedia Authors, 2023)

## Part 4: Conclusions

We have described a universalist view of mathematics and seen that it is in conflict with the literature on the decolonisation of mathematics: the latter is characterised by a questioning of the epistemological privilege given to mathematical knowledge.

Decoloniality and universalism present two very different models of the development of mathematical knowledge. Rather than debate further which is correct, we would suggest that it may be more valuable to apply Box's aphorism and ask which of these models is more useful. The formalist approach we described in the first part of this paper has proved useful, it allowed Gödel to produce profound theorems on the limitations of knowledge, which then led Church and Turing to develop our understanding of the theoretical limits of computation. The rationalist programme of exploring the mathematics of other cultures conducted by Fibonacci and other medieval mathematicians give rise mathematics as we know it today. The programme of ethnomathematics has not, to date, been similarly impactful.

We have also attempted to assemble evidence that might help answer the question of whether the mathematics community is systemically racist. This is of interest whether or not one supports decolonialist approaches. We have looked both decolonialist arguments and empirical evidence. We have focussed on evidence of racism within the UK.

In our view, the theory that the UK curriculum content of mathematics is chosen in a manner that undervalues the contributions of non-European mathematicians is not supported by the evidence. Mathematics is a progressive discipline and the curriculum appears to us to be designed to lead students progressively through the material, and hence it is in large part chronological. The GCSE mathematics syllabus focusses on the contributions of non-European and non-Christian mathematicians, though does not give them explicit credit. We have not checked, but it seems likely that calculus courses will indeed focus on western European mathematicians, but once one moves into the twentieth century the syllabus becomes increasingly diverse.

The evidence put forward to argue that mathematics and colonialism are inseparable is weak, and appears to us primarily rhetorical. Connections can be made between mathematics and colonialism, but given the ubiquity of mathematics in human reasoning this seems unremarkable.

There is some evidence that the contribution of non-Western mathematicians is not sufficiently appreciated by the mathematics community. Decolonialists frequently refer to science and mathematics as western. While authors such as Bishop (1990) should (and indeed do) know better, it does seem plausible that a myth of mathematics as European and colonialist may have gained currency. We can assume that the general public is largely ignorant of mathematical history as the GCSE syllabus does not attempt to teach it. As a result, the public are likely to be readily persuaded by misrepresentation and, unless they consciously choose to study the history of mathematics, they may create their own mythical history of mathematics based on their own biases. However, we did not find evidence that the GCSE syllabus significantly misrepresents the history of mathematics itself.

For these reasons, there is an argument for teaching some aspects of the history of mathematics in schools, or the history of science more generally. However, it is not clear whether that history is best taught by mathematicians or historians. One must also consider the opportunity cost: our informal survey suggests that many in the mathematics faculty at King's College London know very little history of mathematics, yet they function at the highest level as mathematicians. Certainly, if we do teach the history of mathematics, we should endeavour to teach a balanced view of history rather than one designed to support a particular ideological agenda.

The schoolteachers we interviewed were in many ways enthusiastic about teaching a little mathematical history, but they were also very conscious of the importance of teaching the mathematics curriculum itself. As one interviewee pointed out, while the opportunity cost is clear, it is implausible that teaching the history of mathematics would lead to significant benefits for minority students. This view appears to be supported by the evidence on role model interventions (Gladstone & Cimpian, 2021).

The descriptive statistics we have generated on the representation of ethnic minority mathematicians in higher education suggest that there may well be discrimination in science employment in higher education. Perhaps even more compelling than the statistics themselves is the fact that we could not readily find statistics on racial equality in mathematics, and that the London Mathematical Society diversity committee appears to have focussed almost exclusively on gender equality before 2020. Perhaps before we embrace postmodernism, the mathematical community might pause to see what can be gained from a more empirical approach.

## Acknowledgements

Jackman's contribution to this paper was funded by a King's Undergraduate Research Fellowship KURF R4378254. We wish to thank Alan Sokal for his helpful comments on the manuscript.

## Ethical Review

The schoolteacher interviews were approved by the KCL Research Ethics Committee reference number MRSU-22/23-38423.